\documentclass[12pt]{amsart}

\usepackage{amsmath}
\usepackage{amsthm}
\usepackage{amssymb}

\usepackage{amsfonts}

\usepackage[pdftitle={A Foata bijection for the alternating group and for
q analogues},pdfauthor={Dan Bernstein and Amitai Regev},bookmarks=true]{hyperref}

\usepackage{color}

\newcommand{\Z}{\mathbb{Z}}
\newcommand{\st}{\mid} 

\newcommand{\rc}{\mathbf{rc}} %
\newcommand{\rev}{\mathbf{r}} %
\newcommand{\com}{\mathbf{c}} %
\newcommand{\rtlPhi}{\overleftarrow\Phi}    %
\DeclareMathOperator{\Des}{Des} %
\DeclareMathOperator{\Del}{Del} %
\DeclareMathOperator{\del}{del} %
\DeclareMathOperator{\inv}{inv} %
\DeclareMathOperator{\maj}{maj} %
\DeclareMathOperator{\rmaj}{rmaj} %
\DeclareMathOperator{\Pat}{\mathit{Pat}}
\DeclareMathOperator{\Avoid}{\mathit{Avoid}} %
\DeclareMathOperator{\amin}{amin} %
\DeclareMathOperator{\rtlm}{\overleftarrow{\min}} %
\DeclareMathOperator{\ltrm}{\overrightarrow{\min}} %
\DeclareMathOperator{\ltram}{\overrightarrow{\amin}} %

\newtheorem{thm}{Theorem}[section]
\newtheorem{prop}[thm]{Proposition}
\newtheorem{lem}[thm]{Lemma}
\newtheorem{cor}[thm]{Corollary}

\theoremstyle{definition}
\newtheorem{defn}[thm]{Definition}
\newtheorem{rem}[thm]{Remark}
\newtheorem{exmp}[thm]{Example}
\newtheorem{algo}[thm]{Algorithm}

\title{A {F}oata bijection for the alternating group and for
$q$~analogues}

\author{Dan Bernstein}
\address{
    Department of Mathematics\\
    The Weizmann Institute of Science\\
    Rehovot 76100, Israel
}
\email{dan.bernstein@weizmann.ac.il}
\author{Amitai Regev}
\address{
    Department of Mathematics\\
    The Weizmann Institute of Science\\
    Rehovot 76100, Israel
}
\email{am.regev@weizmann.ac.il}

\begin{document}

\date{March 5, 2005}

\begin{abstract}
The Foata bijection $\Phi :S_n\to S_n$ is extended to the
bijections $\Psi :A_{n+1}\to A_{n+1}$ and $\Psi_q :S_{n+q-1}\to
S_{n+q-1}$, where $S_m$, $A_m$ are the symmetric and the
alternating groups. These bijections imply bijective proofs for
recent equidistribution theorems, by Regev and Roichman, for $
A_{n+1}$ and for $S_{n+q-1}$.
\end{abstract}

\maketitle

\section{Introduction}
In~\cite{macmahon16} MacMahon  proved his remarkable theorem
about the equidistribution in the symmetric group $S_n$ of the
{\em length\/} (or the {\em inversion-number\/}) and the {\em major
index\/} statistics. This raised the natural question of
constructing a canonical bijection on $S_n$, for each $n$,  that
would correspond these {\em length\/} and {\em major index\/}
statistics, and would thus yield a bijective proof of that theorem
of MacMahon. That problem was solved by Foata~\cite{foata:netto} ---
who constructed such a canonical bijection, see Section~\ref{SEC:foata} for a discussion of the Foata bijection $\Phi$.
Throughout the years, MacMahon's equidistribution theorem has
received far reaching refinements and generalizations, see for
example~\cite{bjo}, \cite{car1}, \cite{car2}, \cite{foasch},
\cite{gg}, \cite{krat} and~\cite{stanley}.

\medskip

We remark that MacMahon's equidistribution theorem fails when the
$S_n$ statistics are restricted to the alternating subgroups
$A_n$. However, by introducing new $A_n$ statistics which are
natural analogues of the $S_n$ statistics,
in~\cite{regev:alternating}, analogue equidistribution theorems
were proved for $A_n$.  This was done by first formulating the
above $S_n$ statistics in terms of the Coxeter ganerators
$s_i=(i,i+1),\quad 1\le i\le n-1$. By choosing the ``Mitsuhashi''
generators $a_i=s_1s_{i+1}$ for $A_{n+1}$~\cite{mitsuhashi:alternating}, analogue statistics on $A_{n+1}$
were obtained --- via canonical presentations by these generators.
These canonical presentations allow the introduction of the map
$f:A_{n+1}\to S_n$, which is one of the main tools in
~\cite{regev:alternating}, see Section~\ref{pre1} for a discussion
of these presentations and of $f$.

\medskip

MacMahon-type theorems were obtained in~\cite{regev:alternating}
by introducing the {\em delent\/} statistics for these groups (for
$S_n$ --- this statistic already appeared in~\cite{bjo}). Via the
above map $f:A_{n+1}\to S_n$, equidistribution theorems were then
lifted from $S_n$ to $A_{n+1}$, thus yielding (new)
equidistribution theorems for $A_{n+1}$. In particular,
equidistribution theorems for the $A_{n+1}$-analogues of the {\em
length\/} and {\em (reverse) major-index\/} statistics --- were obtained
in this way, see for example Theorem~\ref{a_eq} below.

\medskip

These theorems naturally raise the question of constructing an
$A_{n+1}$-analogue of the Foata bijection --- with the analogue
properties. This problem is solved in Section~\ref{psi}, where
indeed we construct such a map $\Psi$. That map $\Psi$ is composed
of a reflection of Foata's original bijection $\Phi$, together
with the map $f:A_{n+1}\to S_n$ and of certain ``local'' inversions
of $f$. This of course gives a new --- bijective --- proof of Theorem~\ref{a_eq}.

\medskip

Statistics on symmetric groups which are $q$~analogues of the
classical $S_n$ statistics --- were introduced in~\cite{rr2} --- via
$S_n$ canonical presentations and the maps $f_q:S_{n+q-1}\to S_n$,
see Section~\ref{SEC:q_an} below. This map $f_q$ sends the
$q$~statistics on $S_{n+q-1}$ to the corresponding classical
statistics on $S_n$. As in the case of $A_{n+1}$, this allows the
lifting of equidistribution theorems from $S_n$ to $S_{n+q-1}$.
This was done in~\cite{rr2}, where in that process, an interesting
connection with dashed patterns in permutations has appeared, see
Theorems~\ref{qst1} and~\ref{qst2} below.

\medskip

Again, these equidistribution theorems naturally raise the question
of finding the ($q$) analogue of the Foata bijection. In Section~\ref{SEC:q_an} we indeed construct such a bijection $\Psi_q$. As
in the case of $\Psi$, the bijection $\Psi_q$ is composed of
$f_q$, of a reflection of the original Foata bijection $\Phi$, and
of certain ``local'' inversions of $f_q$. As an application,
$\Psi_q$ yields new --- bijective --- proofs of Theorems~\ref{qst1}
and~\ref{qst2}.

\medskip

The paper is organized as follows. Sections~\ref{pre1}, \ref{SEC:foata} and~\ref{pre3} contain
preliminary material, mostly from~\cite{foata:netto},
\cite{regev:alternating} and~\cite{rr2}, which is necessary for
defining and studying the bijections $\Psi$ and $\Psi_q$. A reader
who is familiar with these three papers can skip these preliminary
sections. In Section~\ref{psi} we introduce and study $\Psi$, and Section~\ref{SEC:example}
is an example, showing the properties of $\Psi$. Finally, the
bijections $\Psi_q$ are introduced and studied in Section~\ref{SEC:q_an}.

\section{Canonical presentations and the covering map}\label{pre1}

\subsection{The $S$- and $A$-canonical presentations}\label{SEC:presentations}
In this subsection we
review the presentations of elements in $S_n$ and $A_n$ by the
corresponding generators and procedures for calculating them.

The Coxeter generators of $S_n$ are $s_i=(i,i+1)$, $1\le i\le
n-1$. Recall the definition of the set $R^S_j$,
\[
    R^S_j = \{ 1,\,s_j,\, s_j s_{j-1},\,\dots,\,s_j s_{j-1} \cdots s_1 \} \subseteq S_{j+1} ,
\]
and the following theorem.
\begin{thm}[{see~\cite[pp.~61--62]{goldschmidt:characters}}]
Let $w \in S_n$. Then there exist unique elements $w_j \in R_j^S$, $1 \le j \le n-1$, such that $w = w_1 \cdots w_{n-1}$. Thus, the presentation $w=w_1\cdots w_{n-1}$ is unique. Call that presentation {\em the $S$-canonical presentation of $w$}.
\end{thm}

\medskip

The number of $s_i$ in the $S$-canonical presentation of $\sigma
\in S_n$ is its {\em $S$-length\/}, $\ell_S (\sigma)$. The {\em descent
set\/} of $\sigma\in S_n$ is $\Des_S(\sigma)= \{\,i\mid \sigma
(i)>\sigma (i+1)\,\}$; the  {\em $S$ major-index\/} is
$\maj_S(\sigma)=\sum_{i\in \Des_S(\sigma)}i$, and the {\em reverse $S_n$
major-index\/} is $\rmaj_{S_n}(\sigma)=\sum_{i\in
\Des_S(\sigma)}(n-i)$. Note that $i\in \Des_S(\sigma)$ iff
$\ell_S(\sigma)>\ell_S(\sigma s_i)$.

\medskip

{\bf The $S$-procedure} is a simple way to calculate the
$S$-canonical presentation of a given element in $S_n$. Let
$\sigma \in S_n$, $\sigma(r)=n$, $\sigma = [\dots,n,\dots]$. Then
by definition of the $s_i$s, $n$ can be `pulled to its place on
the right': $\sigma s_r s_{r+1}\cdots s_{n-1} = [\dots, n]$. This
gives $w_{n-1} = s_{n-1}\cdots s_{r+1} s_r \in R^S_{n-1}$. Looking
at $\sigma w_{n-1}^{-1} = \sigma s_r s_{r+1}\cdots s_{n-1} =
[\dots,n-1,\dots,n]$ now, pull $n-1$ to its right place (second
from the right) by a similar product $s_t s_{t+1}\cdots s_{n-2}$,
yielding $w_{n-2} = s_{n-2}\cdots s_t \in R^S_{n-2}$. Continue
this way until finally $\sigma = w_1 \cdots w_{n-1}$.

\begin{exmp}\label{EXMP:s-proc}
Let $\sigma=[5,6,3,2,1,4]$, then $w_5 = s_5 s_4 s_3 s_2$; $\sigma w_5^{-1} = [5,3,2,1,4,6]$, so in order to `pull 5 to its place' we need $w_4 = s_4 s_3 s_2 s_1$; now $\sigma w_5^{-1} w_4^{-1} = [3,2,1,4,5,6]$, so no need to move 4, hence $w_3 = 1$; continuing the same way, check that $w_2 = s_2 s_1$ and $w_1 = s_1$, so $\sigma=w_1 w_2 w_3 w_4 w_5 = (s_5 s_4 s_3 s_2)(s_4 s_3 s_2 s_1)(1)(s_2 s_1)(s_1)$. Thus $\ell_S (\sigma)=11$. Here
$\Des_S\sigma=\{2,3,4\}$, so $\maj_S(\sigma)=\rmaj_{S_6}(\sigma)=9$.
\end{exmp}

\medskip

For $A_n$, the ``Mitsuhash'' generators are $a_i=s_1s_{i+1}$, $1\le i\le n-2$. Recall the definition
\[
    R_j^A = \{1,\,a_j,\,a_j a_{j-1},\,\dots,\,a_j \cdots a_2,\,a_j \cdots a_2 a_1,\,a_j \cdots a_2 a_1^{-1}\} \subseteq A_{j+2}
\]
(for example, $R_3^A = \{1, a_3, a_3 a_2, a_3 a_2 a_1, a_3 a_2 a_1^{-1}\}$), and the following theorem.
\begin{thm}[{see~\cite[Theorem~3.4]{regev:alternating}}]
Let $v \in A_{n+1}$. Then there exist unique elements $v_j \in R_j^A$, $1 \le j \le n-1$, such that $v = v_1 \cdots v_{n-1}$, and this presentation is unique.
Call that presentation {\em the $A$-canonical presentation of $v$}.
\end{thm}

\medskip

The number of $a_i$ in the $A$-canonical presentation of
$\sigma\in A_{n+1}$ is defined to be its {\em $A$-length\/}, $\ell_A (\sigma)$.
In analogy with $S_n$, the $A$-descent set of $\sigma\in A_{n+1}$ is
defined as $\Des_A(\sigma)=\{\,i\mid \ell_A(\sigma)\ge \ell_A(\sigma
a_i)\,\}$. Now define $\maj_A(\sigma)=\sum_{i\in \Des_A(\sigma)}i$,
and $\rmaj_{A_{n+1}}(\sigma)=\sum_{i\in \Des_A(\sigma)}(n-i)$, see~\cite{rr2}.

\medskip

{\bf The $A$-procedure} is a simple procedure for obtaining the
$A$-canonical presentation of $\sigma\in A_n$.\\{\em
$A$-procedure\/}: Step 1: follow the {\em $S$-procedure\/} and
obtain the $S$-canonical presentation of $\sigma\in A_n$. Step 2:
pair the factors. Step 3: insert $s_1s_1$ in the middle of each
pair and obtain the $A$-canonical presentation.

\begin{exmp}\label{EXMP:a-proc}. Let $\sigma = [6,4,3,7,5,2,1]$.\\
Step 1: $\sigma = s_1s_2s_1s_3s_2s_1s_4s_3s_5s_4s_3s_2s_1s_6s_5s_4$.\\
Step 2: $\sigma =
(s_1s_2)(s_1s_3)(s_2s_1)(s_4s_3)(s_5s_4)(s_3s_2)(s_1s_6)(s_5s_4)$.\\
Step 3: $\sigma =
(s_1s_1s_1s_2)(s_1s_1s_1s_3)(s_2s_1s_1s_1)(s_4s_1s_1s_3)(s_5s_1s_1s_4)(s_3s_1s_1s_2)
(s_1s_1s_1s_6)(s_5s_1s_1s_4)=$\\
$=(s_1s_2)(s_1s_3)(s_2s_1)(s_1s_4s_1s_3)(s_1s_5s_1s_4)(s_1s_3s_1s_2)
(s_1s_6)(s_1s_5s_1s_4)=$\\
$=a_1a_2 a_1^{-1}a_3 a_2a_4 a_3 a_2 a_1a_5 a_4 a_3 = (a_1)(a_2
a_1^{-1})(a_3 a_2)(a_4 a_3 a_2 a_1)(a_5 a_4 a_3)$.
\\
Thus $\ell_A(\sigma)=12$ (while $\ell_S(\sigma)=16$). It can be
shown here that $\Des_A(\sigma)=\{1,3,4,5\}$,
hence $\rmaj_{A_7}(\sigma)= 10$.
\end{exmp}

\subsection{The covering map $f$}
We can now introduce the {\em covering map\/} $f$, which plays an
important role in later sections in the constructions of the
bijections $\Psi$ and $\Psi_q$.

\begin{defn}[{see~\cite[Definition~5.1]{regev:alternating}}]
Define $f:R_j^A \to R_j^S$ by
\begin{enumerate}
\item
 $f(a_j a_{j-1}\cdots a_\ell) = s_j s_{j-1}\cdots s_\ell$ if $\ell\ge
 2$, and
 \item
   $f(a_j\cdots a_1) = f(a_j\cdots a_1^{-1}) = s_j\cdots s_1$.
\end{enumerate}
Now extend $f:A_{n+1} \to S_n$ as follows:
let $v\in A_{n+1}$, $v=v_1 \cdots v_{n-1}$ its $A$-canonical presentation, then
\[
    f(v) := f(v_1)\cdots f(v_{n-1}),
\]
which is clearly the $S$-canonical presentation of $f(v)$.
\end{defn}

\section{The Foata bijection}\label{SEC:foata}

The {\em second fundamental transformation on words\/} $\Phi$ was introduced in~\cite{foata:netto} (for a full description, see~\cite[\S10.6]{lothaire:words}). It is defined on any finite word $r=x_1 x_2 \dots x_m$ whose letters $x_1,\dots,x_m$ belong to a totally ordered alphabet.

\begin{defn}\label{DEF:gamma}
Let $X$ be a totally ordered alphabet, let $r=x_1 x_2 \dots x_m$ be a word whose letters belong to $X$, and let $x \in X$ such that $x_m \le x$ (respectively $x_m>x$). Let
\[
    r = r^1 r^2 \dots r^p
\]
be the unique decomposition of $r$ into subwords $r^i=r^i_1 r^i_2
\dots r^i_{m_i}$, $1 \le i \le p$, such that $r^i_{m_i} \le x$
(respectively $r^i_{m_i} > x$) and $r^i_j > x$ (respectively
$r^i_j \le x$) for all $1 \le j < m_i$. Define $\gamma_x(r)$ by
\[
    \gamma_x(r) = r^1_{m_1} r^1_1 r^1_2 \dots r^1_{m_1-1} r^2_{m_2} r^2_1 \dots r^2_{m_2-1}\dots r^p_{m_p} r^p_1 \dots r^p_{m_p-1}.
\]
\end{defn}
For example, with the usual order on the integers, $r=1\,2\,6\,7\,8\,3\,4$ and $x=5$, $r$ decomposes into $r^1=1$, $r^2=2$, $r^3=6\,7\,8\,3$ and $r^4=4$, so
\[
    \gamma_5(1\,2\,6\,7\,8\,3\,4) = 1\,2\,3\,6\,7\,8\,4 .
\]

\begin{defn}
Define $\Phi$ recursively as follows. First, $\Phi(r):=r$ if $r$ is of length 1. If $x$ is a letter and $r$ is a nonempty word, define $\Phi(r x) = \gamma_x(\Phi(r))\,x$.
\end{defn}

For example,
\[
\begin{split}
    \Phi(6\,5\,3\,1\,4\,2)
            &= \gamma_2(\gamma_4(\gamma_1(\gamma_3(\gamma_5(6)\,5)\,3)\,1)\,4)\,2   \\
            &= \gamma_2(\gamma_4(\gamma_1(\gamma_3(6\,5)\,3)\,1)\,4)\,2 \\
            &= \gamma_2(\gamma_4(\gamma_1(6\,5\,3)\,1)\,4)\,2   \\
            &= \gamma_2(\gamma_4(6\,5\,3\,1)\,4)\,2 \\
            &= \gamma_2(3\,6\,5\,1\,4)\,2   \\
            &= 3\,6\,5\,4\,1\,2 .
\end{split}
\]

The following algorithmic description of $\Phi$ from~\cite{foasch}
is more useful in calculations.

\begin{algo}[$\Phi$]\label{ALGO:Phi}
Let $r=x_1 x_2 \dots x_m$ ;

1. Let $i:=1$, $r'_i := x_1$ ;

2. If $i=m$, let $\Phi(r):=r'_i$ and stop; else continue;

3. If the last letter of $r'_i$ is less than or equal to (respectively greater than) $x_{i+1}$, cut $r'_i$ after every letter less than or equal to (respectively greater than) $x_{i+1}$ ;

4. In each compartment of $r'_i$ determined by the previous cuts,
move the last letter in the compartment to the beginning of it;
let $t'_i$ be the word obtained after all those moves; put
$r'_{i+1} := t'_i \, x_{i+1}$ ; replace $i$ by $i+1$ and go to
step 2.
\end{algo}

\begin{exmp}
Calculating $\Phi(r)$, where $r=6\,5\,3\,1\,4\,2$, using the algorithm:
\[
\begin{aligned}
r'_1 &=         6 \mid\\
r'_2 &=         6 \mid 5 \mid\\
r'_3 &=         6 \mid 5 \mid 3 \mid\\
r'_4 &=         6 \;\;\;5\;\;\;3 \mid 1 \mid\\
r'_5 &=         3 \mid 6 \mid 5 \mid 1\;\;\;4 \mid\\
\Phi(r) = r'_6 &=         3\;\;\;6\;\;\;5\;\;\;4\;\;\;1\;\;\;2 \quad.
\end{aligned}
\]
\end{exmp}

The main property of $\Phi$ is the following theorem.

\begin{thm}[see~\cite{foata:netto}]
\label{THM:foata1}
\begin{enumerate}
\item
$\Phi$ is a bijection of $S_n$ onto itself.
\item
For every $\sigma \in S_n$, $\maj_S(\sigma)=\ell_S(\Phi(\sigma))$.
\end{enumerate}
\end{thm}
Some further properties of $\Phi$ are given in Theorem
\ref{THM:foata}

\medskip

Let $\sigma=[\sigma_1,\sigma_2,\dots,\sigma_n] \in S_n$. Denote
the {\em reverse\/} and the {\em complement\/} of $\sigma$ by
\[
    \rev(\sigma):=[\sigma_n,\sigma_{n-1},\dots,\sigma_1]\in S_n
\]
and
\[
    \com(\sigma):=[n+1-\sigma_1,n+1-\sigma_2,\dots,n+1-\sigma_n]\in S_n
\]
respectively.

\begin{rem}\label{REM:rc}
Let $\rho=[n,n-1,\dots,1]\in S_n$. Then for $\sigma \in S_n$,
$\rev(\sigma)=\sigma \rho$ and $\com(\sigma)=\rho \sigma$. Thus it
is obvious that $\rev$ and $\com$ are involutions and that
$\rev\com = \com\rev$. Moreover, $(\rev (\sigma))^{-1} =
\com(\sigma^{-1})$.
\end{rem}

\begin{defn}Let $\rtlPhi := \rev \Phi \rev$, the
{\em right-to-left Foata transformation}.
\end{defn}

While $\rtlPhi(w)$ is easy enough to calculate by reversing $w$,
applying Algorithm~\ref{ALGO:Phi} and reversing the result, it is
easy to see that it can be calculated ``directly'' by applying a
``right-to-left'' version of the Algorithm, namely:

\begin{algo}[$\rtlPhi$]\label{ALGO:rtlPhi}
Let $r=x_1 x_2 \dots x_m$ ;

1. Let $i:=m$, $r'_i := x_m$ ;

2. If $i=1$, let $\Phi(r):=r'_i$ and stop; else continue;

3. If the first letter of $r'_i$ is less than or equal to
(respectively greater than) $x_{m-i}$, cut $r'_i$ before every
letter less than or equal to (respectively greater than) $x_{m-i}$
;

4. In each compartment of $r'_i$ determined by the previous cuts,
move the first letter in the compartment to the end of it; let
$t'_i$ be the word obtained after all those moves; put $r'_{i-1}
:= x_{m-i } \, t'_i$ ; replace $i$ by $i-1$ and go to step 2.
\end{algo}
For an example of applying Algorithm~\ref{ALGO:rtlPhi}, see the
calculation of $\rtlPhi(w)$ in Section~\ref{SEC:example}.

\medskip

The key property of $\rtlPhi$ used in this paper is the following.

\begin{thm}\label{PR:revPhi2}
For every $\sigma \in S_n$, $\rmaj_{S_n}(\sigma) =
\ell_S(\rtlPhi(\sigma))$.
\end{thm}

The proof requires the following lemmas.

\begin{lem}\label{LEM:comPhi}
The bijections $\Phi:S_n\to S_n$ and $\com:S_n\to S_n$ commute with each other.
\end{lem}

\begin{proof}
We prove a slightly stronger claim, namely that $\Phi$ and
$\com_k$ commute as maps on $\Z^n$, where $\com_k(a_1, a_2,\dots,
a_n):=(k+1-a_1,k+1-a_2,\dots,k+1-a_n)$. Let $\sigma=[\sigma_1
\sigma_2 \dots \sigma_n] \in \Z^n$. We need to show that
$\Phi\com_k(\sigma) = \com_k\Phi(\sigma)$. The proof is by
induction on $n$. For $n=1$, everything is trivial. For $n\ge 2$,
write $\Phi(\sigma_1\,\dots\,\sigma_{n-1}) = b_1\,\dots\,b_{n-1}$
and $\gamma_{\sigma_n}(b_1\,\dots\,b_{n-1}) =
c_1\,\dots\,c_{n-1}$. Using the notation $\overline a = k+1-a$, we
have
\[
\begin{split}
\Phi\com_k(\sigma)
    &=\gamma_{\overline{\sigma_n}}( \Phi( \overline{\sigma_1}\, \overline{\sigma_2}\, \dots\, \overline{\sigma_{n-1}}))\,\overline{\sigma_n} \\
    &=\gamma_{\overline{\sigma_n}}( \overline{b_1}\, \overline{b_2}\, \dots\, \overline{b_{n-1}})\,\overline{\sigma_n\vphantom{b}}
\end{split}
\]
by the induction hypothesis, and
\[
\com_k\Phi(\sigma) = \com_k(c_1\,\dots\,c_{n-1}\,\sigma_n) =
\overline{c_1}\,\dots\,\overline{c_{n-1}}\,\overline{\sigma_n} ,
\]
so it remains to show that $\gamma_{\overline{\sigma_n}}(\overline{b_1}\,\overline{b_2}\,\dots\,\overline{b_{n-1}}) = \overline{c_1\vphantom{b}}\,\dots\,\overline{c_{n-1}\vphantom{b}}$.

Assume for now that $b_{n-1}<\sigma_n$ (the case $b_{n-1}>\sigma_n$ is entirely symmetric and will be left to the reader). Let $M=\{\, 1 \le m \le n-1 \st b_m < \sigma_n \,\} = \{m_1,\dots,m_p\}$, $m_1<\dots<m_p$. Note that $\overline{b_{n-1}}>\overline{\sigma_n\vphantom{b_1}}$ and $M=\{\, 1 \le m \le n-1 \st \overline{\sigma_m\vphantom{b_1}} > \overline{b_n} \}$. Therefore, using the notation from Definition~\ref{DEF:gamma}, we have the decompositions
\[
b_1\,\dots\,b_{n-1} = r^1_1\dots r^1_{m_1}\, r^2_1 \dots r^2_{m_2} \dots r^p_1 \dots r^p_{m_p}
\]
and
\[
\overline{b_1\vphantom{^{1}}}\,\dots\,\overline{b_{n-1}\vphantom{^{1}}} =
\overline{r^1_1}\dots \overline{r^1_{m_1\vphantom{1}}}\, \overline{r^2_1} \dots \overline{r^2_{m_2\vphantom{1}}} \dots \overline{r^p_1} \dots \overline{r^p_{m_p\vphantom{1}}}  ,
\]
so
\begin{align*}
    c_1 \dots c_{n-1} = \gamma_{\sigma_n}(b_1\dots b_{n-1})
        &= r^1_{m_1}\, r^1_1 \dots r^1_{m_1-1}\, r^2_{m_2}\, r^2_1 \dots r^2_{m_2-1} \dots r^p_{m_p}\, r^p_1 \dots r^p_{m_p-1}\\
\intertext{and}
    \gamma_{\overline{\sigma_n}}(\overline{b_1}\dots\overline{b_{n-1}}) &= \overline{r^1_{m_1\vphantom{1}}}\, \overline{r^1_1} \dots \overline{r^1_{m_1-1}}\, \overline{r^2_{m_2\vphantom{1}}}\, \overline{r^2_1} \dots \overline{r^2_{m_2-1}} \dots \overline{r^p_{m_p\vphantom{1}}}\, \overline{r^p_1} \dots \overline{r^p_{m_p-1}}   .
\end{align*}
Thus $\gamma_{\overline{\sigma_n}}(\overline{b_1}\,\overline{b_2}\,\dots\,\overline{b_{n-1}}) = \overline{c_1\vphantom{b}}\,\dots\,\overline{c_{n-1}\vphantom{b}}$ as desired.
\end{proof}

\begin{lem}\label{LEM:length}
For every $w \in S_n$, $\ell_S(\rev\com(w)) = \ell_S(w)$.
\end{lem}
\begin{proof}
The lemma follows from the definitions of $\rev$ and $\com$ and from the fact that for all $\sigma \in S_n$, $\ell_S(\sigma) = \inv(\sigma) = \# \{\,(i,j) \st 1\le i<j \le n,\, \sigma(i)>\sigma(j) \,\}$:
\[
\begin{split}
    \ell_S(w) &= \# \{\,(i,j) \st 1\le i<j \le n,\, w(i)>w(j) \,\} \\
              &= \# \{\,(i,j) \st 1\le i<j \le n,\, \com(w)(i)<\com(w)(j) \,\} \\
              &= \# \{\,(i,j) \st 1\le i<j \le n,\, \rev\com(w)(n+1-i)<\rev\com(w)(n+1-j) \,\} \\
              &= \# \{\,(n+1-s,n+1-r) \st 1\le r<s \le n,\, \rev\com(w)(s)<\rev\com(w)(r) \,\} \\
              &= \# \{\,(r,s) \st 1\le r<s \le n,\, \rev\com(w)(r)>\rev\com(w)(s) \,\} \\
              &= \ell_S(\rc(w))      \qedhere
\end{split}
\]
\end{proof}

\begin{lem}\label{LEM:rmaj}
For every $w \in S_n$, $\maj_S(\rev\com(w)) = \rmaj_{S_n}(w)$.
\end{lem}
\begin{proof}
By the definitions of $\rev$, $\com$ and $\Des_S$,
\[
\begin{split}
    i \in \Des_S(\rev\com(w))
        &\iff \rev\com(w)(i)>\rev\com(w)(i+1) \\
        &\iff \com(w)(n-i+1)>\com(w)(n-i) \\
        &\iff n+1-w(n-i+1)>n+1-(w)(n-i) \\
        &\iff w(n-i+1)<(w)(n-i) \\
        &\iff n-i \in \Des_S(w) .
\end{split}
\]
Therefore
\[
    \maj_S(\rev\com(w)) = \sum_{i \in \Des_S(\rev\com(w))} i = \sum_{i \in \Des_S(w)} n-i = \rmaj_{S_n}(w)  . \qedhere
\]
\end{proof}

\begin{proof}[Proof of Theorem~\ref{PR:revPhi2}]
\begin{align*}
    \rmaj_{S_n}(\sigma)
        &= \maj_S(\rc(\sigma))          &&\text{(by Lemma~\ref{LEM:rmaj})}\\
        &= \ell_S(\Phi\rc(\sigma))      &&\text{(by Theorem~\ref{THM:foata1})}\\
        &= \ell_S(\rc\Phi\rc(\sigma))   &&\text{(by Lemma~\ref{LEM:length})}\\
        &= \ell_S(\rtlPhi(\sigma))  &&\text{(by Lemma~\ref{LEM:comPhi} and Remark~\ref{REM:rc})}      \qedhere
\end{align*}
\end{proof}

\section{The delent statistics}\label{pre3}

\begin{defn}[{see~\cite[Definition~7.1]{regev:alternating}}]
Let $\sigma \in S_n$. Define $\Del_S(\sigma)$ as
\[
    \Del_S(\sigma) = \{\, 1<j\le n \st \forall i<j\;\; \sigma(i)>\sigma(j) \,\}.
\]
These are the positions of the l.t.r.min, excluding the first
position.
\end{defn}

\begin{defn}
Let $\sigma \in S_n$. Define the {\em left-to-right minima set\/} of
$\sigma$ as
\[
    \ltrm(\sigma) = \sigma\left(\Del_S(\sigma)\cup\{1\}\right) = \{\, \sigma(j) \st 1\le j \le n,\;\forall i<j\;\;\sigma(i)>\sigma(j) \,\}.
\]
These are the actual (letters) l.t.r.min, including the first
letter.
\end{defn}

\begin{exmp}
Let $\sigma = [5,2,3,1,4]$. Then $\Del_S(\sigma) = \{2,4\}$ and
$\ltrm(\sigma) = \{5,2,1\}$.
\end{exmp}

\begin{prop}
For every $\sigma \in S_n$,
$\ltrm(\sigma)=\Del_S(\sigma^{-1})\cup\{1\}$.
\end{prop}
\begin{proof}
Let $k \in \ltrm(\sigma)$. Then $j = \sigma^{-1}(k) \in
\Del_S(\sigma)\cup\{1\}$. Therefore, by negation, for all $1\le
i\le n$, if $\sigma(i)<\sigma(j) = k$ then $i>j = \sigma^{-1}(k)$.
By the change of variables $i' = \sigma(i)$, we get that for all
$1\le i'\le n$, $i'<k$ implies $\sigma^{-1}(i')>\sigma^{-1}(k)$,
so by definition, $k \in \Del_S(\sigma^{-1})\cup\{1\}$. This
proves that $\ltrm(\sigma) \subseteq
\Del_S(\sigma^{-1})\cup\{1\}$.

The reverse containment is obtained by substituting $\sigma^{-1}$
for $\sigma$ and applying $\sigma$ to both sides.
\end{proof}

\begin{defn}[{see~\cite[Definition~7.4]{regev:alternating}}]
Let $\pi \in A_{n+1}$. Define $\Del_A(\pi)$ as
\[
    \Del_A(\pi) = \{\, 2<j\le n+1 \st \text{there is at most one $i<j$ such that $\pi(i)<\pi(j)$} \,\}.
\]
\end{defn}

\begin{defn}
Let $\pi \in A_{n+1}$. Define the {\em left-to-right almost-minima
set\/} of $\pi$ as
\begin{multline*}
    \ltram(\pi) = \pi\left(\Del_A(\pi)\cup\{1,2\}\right)\\
    =\{\,\pi(j) \st 1 \le j \le n+1 \text{ and there is at most one $i<j$ such that $\pi(i)<\pi(j)$} \,\}.
\end{multline*}
\end{defn}

\begin{exmp}
Let $\pi = [4,2,6,3,1,5]$. Then $\Del_A(\pi) = \{4,5\}$ and
$\ltram(\pi) = \{4,2,3,1\}$.
\end{exmp}

\begin{prop}\label{PR:DelA}
For every $\pi \in A_{n+1}$,
$\ltram(\pi)=\Del_A(\pi^{-1})\cup\{1,2\}$.
\end{prop}
\begin{proof}
Let $k \in \ltram(\pi)$. Then $j = \pi^{-1}(k) \in
\Del_A(\pi)\cup\{1,2\}$. Therefore for all $1\le i\le n+1$ except
at most one, if $\pi(i)<\pi(j) = k$ then $i>j = \pi^{-1}(k)$. By
the change of variables $i' = \pi(i)$, we get that for all $1\le
i'\le n+1$ except at most one, $i'<k$ implies
$\pi^{-1}(i')>\pi^{-1}(k)$, so by definition, $k \in
\Del_A(\pi^{-1})\cup\{1,2\}$. This proves that $\ltrm(\pi)
\subseteq \Del_A(\pi^{-1})\cup\{1,2\}$.

The reverse containment is obtained by substituting $\pi^{-1}$ for
$\pi$ and applying $\pi$ to both sides.
\end{proof}

We now quote the following theorem. The bijection $\Psi$ of
Theorem~\ref{PR:Psi} bellow provides a (short) bijective proof of
that theorem.

\begin{thm}[see~{\cite[Theroem~9.1(2)]{regev:alternating}}]
\label{a_eq}
For every subsets $D_1 \subseteq \{1,\dots,n-1\}$ and $D_2
\subseteq \{1,\dots,n-1\}$,
\[
\sum\limits_{\{\sigma\in A_{n+1} \st  \Des_A(\sigma^{-1})\subseteq
D_1,\ \Del_A(\sigma^{-1})\subseteq D_2\}}
q^{\rmaj_{A_{n+1}}(\sigma)} = \sum\limits_{\{\sigma\in A_{n+1} \st
\Des_A(\sigma^{-1})\subseteq D_1,\ \Del_A(\sigma^{-1})\subseteq
D_2\}} q^{\ell_A(\sigma)}.
\]
\end{thm}

\section{The bijection $\Psi$}\label{psi}

Recall the notations for the reverse and the complement of
$\sigma=[\sigma_1 \sigma_2 \dots \sigma_n] \in S_n$, which are
$\rev(\sigma) = [\sigma_n \sigma_{n-1} \dots \sigma_1]$ and
$\com(\sigma) = [n+1-\sigma_1, n+1-\sigma_2,\dots,n+1-\sigma_n]$,
respectively, and the notations $\Phi$ and $\rtlPhi = \rev \Phi
\rev$ for Foata's {\em second fundamental transformation\/} and
the {\em right-to-left Foata transformation\/} (both described in
detail in Section~\ref{SEC:foata}), respectively.

We shall need the following properties of $\Phi$ and $\rtlPhi$,
see also Theorem~\ref{THM:foata1}.
\begin{thm}
\label{THM:foata}
\begin{enumerate}
\item
$\Phi$ is a bijection of $S_n$ onto itself.
\item
For every $\sigma \in S_n$, $\maj_S(\sigma)=\ell_S(\Phi(\sigma))$.
\item
{\rm (see~\cite[Example~5.3]{bjo})} For every $\sigma \in S_n$,
$\rtlm(\sigma)=\rtlm(\Phi(\sigma))$, where $\rtlm(\sigma) = \{\,
\sigma(j) \st 1\le j \le n,\;\forall i>j\;\;\sigma(i)>\sigma(j)
\,\}$ is the set of right-to-left minima of $\sigma$.
\item
{\rm (see~\cite[Theorem~1]{foasch})} For every $\sigma \in S_n$,
$\Des_S(\sigma^{-1})=\Des_S([\Phi(\sigma)]^{-1})$.
\item
By Theorem~\ref{PR:revPhi2}, for every $\sigma \in S_n$,
$\rmaj_{S_n}(\sigma) = \ell_S(\rtlPhi(\sigma))$.
\end{enumerate}
\end{thm}

The $S$- and the $A$-canonical presentations and the map $f$ were
discussed in Section~\ref{pre1}. A key property of $f$ is the way
it relates between certain pairs of statistics on $A_{n+1}$ and on
$S_n$.

\begin{defn}[{see~\cite[Definition~5.2]{rr2}}]
Let $m_S$ be a statistic on the symmetric groups and $m_A$ a statistic on the alternating groups. We say that $(m_S, m_A)$ is an {\em $f$-pair (of statistics)\/} if for any $n$ and $v \in A_{n+1}$, $m_A(v) = m_S(f(v))$.
\end{defn}

\begin{prop}[{see~\cite[Propositions~5.3 and~5.4]{regev:alternating}}]\label{PR:fPairs}
The following pairs are $f$-pairs: $(\ell_S,\ell_A)$,
$(\rmaj_{S_n},\rmaj_{A_{n+1}})$, $(\del_S,\del_A)$ and
$(\Des_A,\Des_S)$.
\end{prop}

We also have
\begin{prop}[{see~\cite[Propositions~8.4 and~8.5]{rr2}}]\label{PR:fInv}
For every $v\in A_{n+1}$, $f(v)^{-1} = f(v^{-1})$.
\end{prop}

The covering map $f$ is obviously not injective. The family of
maps $g_u$ defined next serve as local inverses of $f$ (see Remark~\ref{REM:gInvertible}).

\begin{defn}
For $u \in A_{n+1}$ with $A$-canonical presentation $u= u_1 u_2 \cdots u_{n-1}$,
define $g_u:R_j^S \to R_j^A$ by
\[
    g_u(s_j s_{j-1}\cdots s_\ell) = a_j a_{j-1} \cdots a_\ell
    \quad \text{if \;$\ell\ge 2$,\; and} \quad
        g_u(s_j s_{j-1}\cdots s_1) = u_j.
\]
Now extend $g_u:S_n \to A_{n+1}$ as follows:
let $w\in S_n$, $w=w_1 \cdots w_{n-1}$ its $S$-canonical presentation, then
\[
    g_u(w) := g_u(w_1)\cdots g_u(w_{n-1}),
\]
which is clearly the $A$-canonical presentation of $g_u(w)$.
\end{defn}

\begin{rem}\label{REM:gInvertible}
Let $w \in S_n$ and $u\in A_{n+1}$. Then $f(g_u(w))=w$ if for all $1 \le j \le n-1$,
\[
    u_j = a_j\cdots a_2 a_1^{\pm 1} \; \iff \; w_j = s_j\cdots s_1  ,
\]
where $w=w_1\cdots w_{n-1}$ and $u=u_1\cdots u_{n-1}$ are the $S$- and $A$-canonical presentations of $w$ and $u$ respectively.
\end{rem}

We are now ready to define the bijection $\Psi$.

\begin{defn}
Define $\Psi:A_{n+1} \to A_{n+1}$ by $\Psi(v) =
g_v(\rtlPhi(f(v)))$ .
\end{defn}

That is, the image of $v$ under $\Psi$ is obtained by applying
$\rtlPhi$ to $f(v)$ in $S_n$, then using $g_v$ as an ``inverse''
of $f$ in order to ``lift'' the result back to $A_{n+1}$.

The following is our main theorem, which can be seen as an
$A_{n+1}$-analogue of Theorem~\ref{THM:foata}.

\begin{thm}\label{PR:Psi}
\begin{enumerate}
\item
The mapping $\Psi$ is a bijection of $A_{n+1}$ onto itself.
\item
For every $v\in A_{n+1}$, $\rmaj_{A_{n+1}}(v) = \ell_A(\Psi(v))$.
\item
For every $v\in A_{n+1}$, $\del_A(v)=\del_A(\Psi(v))$.
\item
For every $v\in A_{n+1}$, $\Del_A(v^{-1}) =
\Del_A([\Psi(v)]^{-1})$.
\item
For every $v\in A_{n+1}$, $\Des_A(v^{-1})=\Des_A([\Psi(v)]^{-1})$.
\end{enumerate}
\end{thm}

In order to prove the theorem we need the following lemmas.

\begin{lem}\label{LEM:rep-ltrm}
\begin{enumerate}
\item
Let $w \in S_n$, $w=w_1 \cdots w_{n-1}$ its $S$-canonical
presentation. Then for every $1<j \le n$, $j\in\ltrm(w)$ if and
only if $w_{j-1} = s_{j-1} s_{j-2}\cdots s_1$.
\item
Let $v \in A_{n+1}$, $v=v_1 \cdots
v_{n-1}$ its $A$-canonical presentation. Then for every $2<j \le
n+1$, $j\in\ltram(v)$ if and only if $v_{j-2} = a_{j-2}
a_{j-3}\cdots a_1^{\pm 1}$.
\end{enumerate}
\end{lem}
\begin{proof}
\begin{enumerate}
\item
By induction on $n$. Let $\sigma = w_1 \cdots w_{n-2} \in S_{n-1}
\subseteq S_n$ and assume that the assertion is true for $\sigma$.
If $w_{n-1}=1$, then the claim is correct by the induction
hypothesis. Otherwise, $w_{n-1}=s_{n-1}s_{n-2}\cdots s_\ell$ for
some $1 \le \ell \le n-1$. Writing $\sigma=[b_1,\dots,b_{n-1}]$,
we have that $w=\sigma w_{n-1} =
[b_1,\dots,b_{\ell-1},n,b_\ell,\dots,b_{n-1}]$. For every $1<j\le
n-1$, $j=b_k$ for some $k$, so $j$ is a left-to-right minimum of
$w$ if and only if it is a left-to-right minimum of $\sigma$,
which, by the induction hypothesis, is true if and only if
$w_{j-1}=s_{j-1}\cdots s_1$. Finally, $n$ is an additional
left-to-right minimum of $w$ if and only if $\ell=1$, that is if
and only if $w_{n-1}=s_{n-1}s_{n-2}\cdots s_1$.
\item
By induction on $n$. Let $\pi = v_1 \cdots v_{n-2} \in A_{n}
\subseteq A_{n+1}$ and assume that the assertion is true for
$\pi$. If $v_{n-1}=1$, then the claim is correct by the induction
hypothesis. Otherwise, $v_{n-1}=a_{n-1}a_{n-2}\cdots
a_\ell^\epsilon$ for some $1 \le \ell \le n-1$ and $\epsilon=\pm
1$. Writing $\pi=[c_1,c_2,\dots,c_n]$, we have that
\[
    v=\pi v_{n-1} = \begin{cases}
        [c_1 c_2,\dots,c_\ell,n+1,c_{\ell+1},\dots,c_n],
            &\text{if $\ell>1$ and $n-\ell$ is even;} \\
        [c_2 c_1,\dots,c_\ell,n+1,c_{\ell+1},\dots,c_n],
            &\text{if $\ell>1$ and $n-\ell$ is odd;} \\
        [c_1,n+1,c_2,\dots,c_n],
            &\text{if $\ell=1$, $n$ is odd and $\epsilon=1$}; \\
        [c_2,n+1,c_1,\dots,c_n],
            &\text{if $\ell=1$, $n$ is even and $\epsilon=1$}; \\
        [n+1,c_1,c_2,\dots,c_n],
            &\text{if $\ell=1$, $n$ is even and $\epsilon=-1$}; \\
        [n+1,c_2,c_1,c_3,\dots,c_n],
            &\text{if $\ell=1$, $n$ is odd and $\epsilon=-1$}. \\
    \end{cases}
\]
For every $2<j\le n$, $j=c_k$ for some $k$, so $j$ is a
left-to-right almost-minimum of $v$ if and only if it is a
left-to-right almost-minimum of $\pi$, which, by the induction
hypothesis, is true if and only if $v_{j-2}=a_{j-2}\cdots a_1^{\pm
1}$. Finally, $n+1$ is an additional left-to-right almost-minimum
of $v$ if and only if $\ell=1$, that is if and only if
$v_{n-1}=a_{n-1}a_{n-2}\cdots a_1^{\pm 1}$.\qedhere
\end{enumerate}
\end{proof}

\begin{cor}\label{nosaf}
For every $v\in A_{n+1}$, $\ltram(v)=\ltrm(f(v))-1$, where
$X-1=\{x-1 \st x\in X\}$.
\end{cor}

\begin{lem}\label{LEM:ltrm}
For every $w \in S_n$, $\ltrm(w) = \ltrm( \rtlPhi(w) )$, hence $\del_S(w)=\del_S(\rtlPhi(w))$.
\end{lem}
\begin{proof}
This follows immediately from the definitions and from Theorem~\ref{THM:foata}(3):
\[
\begin{split}
    j \in \ltrm(w) & \iff j \in \rtlm(\rev(w)) \\
                    & \iff j \in \rtlm(\Phi(\rev(w)) \\
                    & \iff j \in \ltrm(\rev(\Phi(\rev(w))) = \ltrm(\rtlPhi(w)).  \qedhere
\end{split}
\]
\end{proof}

The following is an easy corollary of Lemmas~\ref{LEM:rep-ltrm}
and~\ref{LEM:ltrm}.

\begin{cor}\label{COR:revPhi}
Let $w\in S_n$, $w=w_1\cdots w_{n-1}$ its $S$-canonical
presentation, and let $\sigma=\rtlPhi(w)$, $\sigma=\sigma_1\cdots
\sigma_{n-1}$ its $S$-canonical presentation. Then
$\sigma_j=s_j\cdots s_1$ if and only if $w_j=s_j\cdots s_1$.
\end{cor}

\begin{lem}\label{LEM:fPsi}
Let $v \in A_{n+1}$. Then $f(\Psi(v)) = \rtlPhi(f(v))$.
\end{lem}
\begin{proof}
Let
$v = v_1\cdots v_{n-1}$ and $w=\rtlPhi(f(v))=w_1\cdots w_{n-1}$ be the $A$- and $S$-canonical presentations of $v$ and $\rtlPhi(f(v))$ respectively. By definition of $f$
and Corollary~\ref{COR:revPhi}, for every $1 \le j \le n-1$, $w_j
= s_j s_{j-1}\cdots s_1$ if and only if $v_j = a_j\cdots a_2
a_1^{\pm 1}$. Therefore, by Remark~\ref{REM:gInvertible},
\[
    f(\Psi(v))=f(g_v(\rtlPhi(f(v))))=f(g_v(w))=w = \rtlPhi(f(v)) . \qedhere
\]
\end{proof}

\begin{proof}[Proof of Theorem~\ref{PR:Psi}]
\begin{enumerate}
\item
To prove that $\Psi$ is a bijection, it suffices to find its inverse.
Let $v \in A_{n+1}$, and let $v=v_1\cdots v_{n-1}$, $w=\rtlPhi(f(v))=w_1\cdots w_{n-1}$ and $u=\Psi(v)=g_v(w)=u_1\cdots u_{n-1}$
be the $A$-, $S$- and $A$-canonical presentations of $v$,
$\rtlPhi(f(v))$ and $\Psi(v)$ respectively.
By Lemma~\ref{LEM:fPsi},
\[
\rtlPhi^{-1}(f(\Psi(v))) = \rtlPhi^{-1}(\rtlPhi(f(v)))=f(v),
\]
so
\begin{equation}\label{EQ:inverse}
    g_{\Psi(v)}(\rtlPhi^{-1}(f(\Psi(v)))) = g_{\Psi(v)}(f(v))
  = g_u(f(v)) = g_u(f(v_1))\cdots g_u(f(v_{n-1}))  \tag{$*$}     .
\end{equation}
We claim that $\pi \mapsto g_\pi(\rtlPhi^{-1}(f(\pi)))$ is the inverse of $\Psi$, or in other words, that the right hand side of~\eqref{EQ:inverse} equals $v_1 v_2 \cdots v_{n-1}$. Let $1 \le j \le n-1$. If $v_j = a_j a_{j-1} \cdots a_\ell$, $\ell>1$, then $g_u(f(v_j))=g_u(s_j s_{j-1}\cdots s_\ell)=a_j a_{j-1} \cdots a_\ell=v_j$. If $v_j = a_j\cdots a_2 a_1^{\pm 1}$, then $f(v_j)=s_j s_{j-1}\cdots s_1$, so by Corollary~\ref{COR:revPhi}, $w_j=s_j s_{j-1}\cdots s_1$, and therefore $u_j=g_v(w_j)=v_j$, so again $g_u(f(v_j))=v_j$, and the claim is proved.

\item
By Proposition~\ref{PR:fPairs} and Lemma~\ref{LEM:fPsi},
$\ell_A(\Psi(v)) = \ell_S(f(\Psi(v))) = \ell_S(\rtlPhi(f(v)))$. By
Theorem~\ref{PR:revPhi2} and Proposition~\ref{PR:fPairs},
$\ell_S(\rtlPhi(f(v))) = \rmaj_{S_n}(f(v)) = \rmaj_{A_{n+1}}(v)$.
Thus $\ell_A(\Psi(v)) = \rmaj_{A_{n+1}}(v)$ as desired.

\item
By Proposition~\ref{PR:fPairs} and Lemma~\ref{LEM:fPsi},
$\del_A(\Psi(v))=\del_S(f(\Psi(v)))=\del_S(\rtlPhi(f(v)))$, and by
Lemma~\ref{LEM:ltrm}, the definition of $\del_S$ and
Proposition~\ref{PR:fPairs}, $\del_S(\rtlPhi(f(v))) = \del_S(f(v))
= \del_A(v)$. Thus $\del_A(\Psi(v))=\del_A(v)$ as desired.
\item

By Corollary~\ref{nosaf},
$\ltram(\Psi(v)) = \ltrm(f(\Psi(v)))-1$, with the notation $X-1 =
\{\,x-1 \st x\in X\,\}$. Therefore by Lemmas~\ref{LEM:fPsi}
and~\ref{LEM:ltrm}, $\ltram(\Psi(v)) = \ltrm(\rtlPhi(f(v)))-1 =
\ltrm(f(v))-1$. Again by Lemma~\ref{LEM:rep-ltrm}, we get that
$\ltram(\Psi(v)) = \ltram(v)$. By Proposition~\ref{PR:DelA}, this
implies that $\Del_A([\Psi(v)]^{-1})\cup\{1,2\} =
\Del_A(v^{-1})\cup\{1,2\}$, hence $\Del_A([\Psi(v)]^{-1}) =
\Del_A(v^{-1})$ as desired.

\item
By Propositions~\ref{PR:fPairs} and~\ref{PR:fInv} and Lemma~\ref{LEM:fPsi},
\[\Des_A([\Psi(v)]^{-1}) = \Des_S(f([\Psi(v)]^{-1})) = \Des_S([f(\Psi(v))]^{-1})
\Des_S([\rtlPhi(f(v))]^{-1}).\] By Remark~\ref{REM:rc},
$\rtlPhi(f(v))^{-1} = (\rev \Phi \rev f(v))^{-1} = \com((\Phi \rev
f(v))^{-1})$, so $\Des_S([\rtlPhi(f(v))]^{-1}) = \{1,\dots,n-1\}
\setminus \Des_S([\Phi \rev f(v)]^{-1})$. By
Theorem~\ref{THM:foata}, \[\Des_S([\Phi\rev f(v)]^{-1}) =
\Des_S([\rev f(v)]^{-1}).\] Hence, $\Des_S([\rtlPhi(f(v))]^{-1}) =
\{1,\dots,n-1\} \setminus \Des_S([\rev f(v)]^{-1}) =
\Des_S(\com([\rev f(v)]^{-1}))$. Since $\com([\rev f(v)]^{-1}) =
\com(\com([f(v)]^{-1})) = f(v)^{-1}$, we get that
\[\Des_S([\rtlPhi(f(v))]^{-1}) = \Des_S([f(v)]^{-1}).\] Finally, by
Propositions~\ref{PR:fInv} and~\ref{PR:fPairs},
$\Des_S([f(v)]^{-1})=\Des_S(f(v^{-1}))=\Des_A(v^{-1})$. \qedhere
\end{enumerate}

\end{proof}

\section{Example}\label{SEC:example}

As an example, let $v=[6, 4, 3, 7, 5, 2, 1] \in A_7$. We now
calculate $v$, $v^{-1}$, $\Psi (v)$ and $[\Psi (v)]^{-1}$, and
using the $A$-procedure --- their $A$-canonical presentations.
This yields the corresponding sets $\Del_A$ and $\Des_A$, hence
also the $\ell_A$ and the $rmaj_{A_7}$ indices, thus demonstrating
Theorem~\ref{PR:Psi} in this example. Throughout the example, when
writing a canonical presentation, we will underline all factors of
the form $\underline{a_j\cdots a_2 a_1^{\pm 1}}$ and
$\underline{s_j \cdots s_1}$.

\medskip

The $A$-canonical presentations of $v$ and of $v^{-1}$ are
\[
    v=\underline{v_1} \underline{v_2} v_3 \underline{v_4} v_5 =
    (\underline{a_1
    \vphantom{a_1^{}}})
    (\underline{a_2 a_1^{-1}})(a_3 a_2)(\underline{a_4 a_3 a_2 a_1\vphantom{a_1^{}}})(a_5 a_4
    a_3)
\quad\text{(so\quad$\del_A(v) = 3$)},\]
\[v^{-1}=[7,6,3,2,5,1,4]=
(\underline{a_1\vphantom{a_1^{}}})(a_3a_2)(\underline{a_4a_3a_2a_1^{-1}})
(\underline{a_5a_4a_3a_2a_1^{-1}})\quad\text{(so\quad$\del_A(v^{-1})=3)$}.\] Thus $\Des_A(v)=\{1,3,4,5\}$, so
$\rmaj_{A_7}(v)=(6-1)+(6-3)+(6-4)+(6-5) = 11$.\\
Similarly $\Des_A(v^{-1})=\{1,2,4\}$. Also, $\Del_A(v)=\{3,6,7\}$
and $\Del_A(v^{-1})=\{3,4,6\}$.

\medskip

 We have
\[
    w = f(v) = \underline{w_1}\underline{w_2}w_3\underline{w_4}w_5 = (\underline{s_1})(\underline{s_2 s_1})(s_3 s_2)(\underline{s_4 s_3 s_2 s_1})(s_5 s_4 s_3) = [5,3,6,4,2,1] .
\]
Note that $\Des_S(w)=\Des_S(f(v))= \{1,3,4,5\}=\Des_A(v)$, and also,
$\rmaj_{S_6}(w) = 11 = \rmaj_{A_7}(v)$ and
$\del_S(w)=3=\del_A(v)$, in accordance with
Proposition~\ref{PR:fPairs}.

\medskip

Let us calculate $\Psi(v)$ and $[\Psi(v)]^{-1}$. Using
Algorithm~\ref{ALGO:rtlPhi} we obtain $\rtlPhi(w)$:
\[
\begin{aligned}
w'_1 &=         & \mid 1\\
w'_2 &=         & \mid 2 \mid 1\\
w'_3 &=         & \mid 4 \mid 2 \mid 1\\
w'_4 &=         & \mid 6 \mid 4\;\;\;2\;\;\;1\\
w'_5 &=         & \mid 3\;\;\;6 \mid 2 \mid 1 \mid 4\\
\rtlPhi(w) = w'_6 &=         & 5,\,\;6,\,\;3,\,\;2,\,\;1,\,\;4 &.
\end{aligned}
\]
Note that $\ell_S(\rtlPhi(w))=11=\rmaj_{S_6}(w)$, as asserted by
Theorem~\ref{PR:revPhi2}.

The $S$-canonical presentation of $\rtlPhi(w)$, obtained by the
$S$-procedure (see Example~\ref{EXMP:s-proc}), is
\[
    u = \rtlPhi(w) = \underline{u_1}\underline{u_2}u_3\underline{u_4}u_5 = (\underline{s_1})(\underline{s_2 s_1})(1)(\underline{s_4 s_3 s_2 s_1})(s_5 s_4 s_3 s_2)   .
\]
The underlined factors in the $S$-canonical presentation of $w$
are the same as the underlined factors in the $S$-canonical
presentation of $\rtlPhi(w)$, as asserted by
Corollary~\ref{COR:revPhi}. This is a result of the fact that
$\ltrm(w)=\{1,2,3,5\}=\ltrm(\rtlPhi(w))$, which is a result
of Lemma~\ref{LEM:ltrm}.

Now
\begin{multline*}
    \Psi(v)=g_v(u) = \underline{v_1}\underline{v_2}(1)\underline{v_4}(a_5 a_4 a_3 a_2) =\\
    (\underline{a_1\vphantom{a_1^{}}})(\underline{a_2 a_1^{-1}})(1)
    (\underline{a_4 a_3 a_2 a_1\vphantom{a_1^{}}})(a_5 a_4 a_3 a_2) = [4,6,7,3,2,1,5]        ,
\end{multline*}
so $[\Psi(v)]^{-1}=[6,5,4,1,7,2,3]=(1)(\underline{a_2a_1\vphantom{a_1^{}}})
(\underline{a_3a_2a_1\vphantom{a_1^{}}})(\underline{a_4a_3a_2a_1^{-1}})(a_5a_4)$.
It follows that \[\Des_A(v^{-1})=\{1,2,4\}=\Des_A([\Psi(v)]^{-1})
\quad\text{and}\quad
\Del_A(v^{-1})=\{3,4,6\}=\Del_A([\Psi(v)]^{-1}).\] Also
\[
    \del_A(\Psi(v)) = 3 = \del_A(v)
\]
and
\[
    \ell_A(\Psi(v)) = 11 = \rmaj_{A_7}(v).
\]

\section{$q$ analogues}\label{SEC:q_an}
\subsection{The $q$ statistics}

\begin{defn}[{see~\cite[Definition~5.1]{rr2}}]
Let $\pi \in S_n$, and let $q<n$. Define the {\em $q$-length\/} of
$\pi$, $\ell_q(\pi)$, as the number of Coxeter generators in the
$S$-canonical presentation of $\pi$, where $s_1,\dots,s_{q-1}$ are
not counted. For example, let
$\pi=s_1s_2s_1s_4s_3s_6s_5s_4s_3s_2$, then $\ell_3 (\pi)=6$ while
$\ell_4(\pi)=4$. Clearly, $\ell_1=\ell_S$.
\end{defn}

\begin{defn}[{see~\cite[Definition~5.1]{rr2}}]
Let $\pi \in S_n$. Define $\Del_{k+1}(\pi)$ as
\[
    \Del_{k+1}(\pi) = \{\, k+1<j\le n \st \#\{i<j \st \pi(i)<\pi(j)\}\le k \,\}.
\]
\end{defn}
\begin{defn}
Let $\pi \in S_n$. Define the {\em left-to-right $k$-almost-minima
set\/} of $\pi$ as
\begin{multline*}
    \ltrm_{k+1}(\pi) = \pi\left(\Del_{k+1}(\pi)\cup\{1,2,\dots,k+1\}\right)\\
    =\{\,\pi(j) \st 1 \le j \le n,\; \#\{i<j \st \pi(i)<\pi(j)\}\le k \,\}.
\end{multline*}
\end{defn}

\begin{prop}\label{PR:Delk}
For every $\pi \in S_{n+q-1}$,
$\ltrm_{k+1}(\pi)=\Del_{k+1}(\pi^{-1})\cup\{1,2,\dots,k+1\}$.
\end{prop}
\begin{proof}
Let $r \in \ltrm_{k+1}(\pi)$. Then $j = \pi^{-1}(r) \in
\Del_{k+1}(\pi)\cup\{1,\dots,k+1\}$. Therefore
\[
    \#\{ 1 \le i \le n+q-1 \st \pi(i)<\pi(j)=r \text{ and } i<j=\pi^{-1}(r) \} \le k    .
\]
By
the change of variables $i' = \pi(i)$, we get
\[
    \#\{ 1 \le i' \le n+q-1 \st i'<r \text{ and } \pi^{-1}(i')<\pi^{-1}(r) \} \le k     ,
\]
so by definition, $r \in
\Del_{k+1}(\pi^{-1})\cup\{1,\dots,k+1\}$. This proves that $\ltrm_{k+1}(\pi)
\subseteq \Del_{k+1}(\pi^{-1})\cup\{1,\dots,k+1\}$.

The reverse containment is obtained by substituting $\pi^{-1}$ for
$\pi$ and applying $\pi$ to both sides.
\end{proof}

\begin{defn}[{see~\cite[Definition~5.8]{rr2}}]
Let $\pi \in S_{n+q-1}$. Then $i$ is a {\em $q$-descent\/} in $\pi$
if $i\ge q$ and at least one of the following holds: a) $i \in \Des(\pi)$; b) $i+1 \in \Del_q(\pi)$.
\end{defn}

\begin{defn}[{see~\cite[Definition~5.9]{rr2}}]
\begin{enumerate}
\item
The {\em $q$-descent set\/} of $\pi \in S_{n+q-1}$ is defined as
\[
    \Des_q(\pi) = \{\,i \st \text{$i$ is a $q$-descent in
    $\pi$}\,\}.
\]
\item
For $\pi \in S_{n+q-1}$ define the {\em $q,m$-reverse major index\/}
of $\pi$ by
\[
    \rmaj_{q,m}(\pi) = \sum_{i\in\Des_q(\pi)} (m-i)    ,
\]
where $m=n+q-1$.
\end{enumerate}
\end{defn}
We need the notion of {\em dashed\/} patterns~\cite{rr2}, and we
introduce it via examples:\\ $\sigma\in S_n$ has the dashed pattern
$(1-2-4,3)$ if $\sigma=[\cdots,a,\cdots,b,\cdots ,d,c,\cdots]$,
and it has the dashed pattern $(2-1-4,3)$ if
$\sigma=[\cdots,b,\cdots,a,\cdots ,d,c,\cdots]$ for some
$a<b<c<d$. Given $q$, denote by $\Pat(q)$ the following $q!$
dashed patterns: \[\Pat(q) =\{\,(\pi_1-\pi_2-\cdots
-\pi_q-(q+2),(q+1))\mid \pi\in S_q\,\}.\] For example, $\Pat(2)
=\{(1-2-4,3),\;(2-1-4,3)\}$. If $\sigma\in S_m$ does not have any
of the dashed pattern in $\Pat(q)$, then $\sigma$  {\em avoids\/}
$\Pat(q)$. We denote by $\Avoid_q(n+q-1)$ the set of permutations
$\sigma\in S_{n+q-1}$ avoiding all the $q!$ dashed patterns in
$\Pat(q)$.

\medskip

The main equidistribution theorems here are the following two
theorems. The bijection $\Psi_q$ below implies bijective proofs
for these theorems.
\begin{thm}[{see~\cite[Theorem~11.5]{rr2}}]\label{qst1}
For every positive integers $n$ and $q$ and every subsets $B_1,
B_2 \subseteq \{q,\dots,n+q-1\}$,
\[
    \sum_{\{\,\pi \in S_{n+q-1} \st
    \Des_q(\pi^{-1})=B_1,\;\Del_q(\pi^{-1})=B_2\,\}} t^{\ell_q(\pi)}
    =    \sum_{\{\,\pi \in S_{n+q-1} \st
    \Des_q(\pi^{-1})=B_1,\;\Del_q(\pi^{-1})=B_2\,\}} t^{\rmaj_{q,n+q-1}(\pi)}
    .
\]
\end{thm}

\begin{thm}[{see~\cite[Theorem~11.7]{rr2}}]\label{qst2}
For every positive integers $n$ and $q$ and every subsets $B
\subseteq \{q,\dots,n+q-2\}$,
\[
    \sum_{\{\,\pi^{-1} \in \Avoid_q(n+q-1) \st
    \Des_q(\pi^{-1})=B\,\}} t^{\ell_q(\pi)}
    =    \sum_{\{\,\pi^{-1} \in \Avoid_q(n+q-1) \st
    \Des_q(\pi^{-1})=B\,\}} t^{\rmaj_{q,n+q-1}(\pi)}
    .
\]
\end{thm}

\subsection{The covering map $f_q$}
\begin{defn}[{see~\cite[Definition~8.1]{rr2}}]
Let $w \in S_{n+q-1}$ and let $w= s_{i_1}\cdots s_{i_r}$ be its
$S$-canonical presentation. Define $f_q:S_{n+q-1}\to S_n$ as
follows:
\[
    f_q(w) = f_q(s_{i_1})\cdots f_q(s_{i_r}),
\]
where $f_q(s_1) = \cdots = f_q(s_{q-1}) = 1$, and
$f_q(s_j)=s_{j-q+1}$ if $j\ge q$.
\end{defn}
\begin{rem}
If $w = w_1 \cdots w_{n+q-2}$ is the $S$-canonical presentation of
$w\in S_{n+q-1}$, $w_j \in R^S_j$, then $f_q(w) = f_q(w_q)\cdots
f_q(w_{n+q-2})$ is the $S$-canonical presentation of $f_q(w)$,
$f_q(w_j) \in R^S_{j-q+1}$.
\end{rem}

\begin{prop}[{see~\cite[Proposition~8.6 and Remark~11.1]{rr2}}]\label{PR:fqPairs}
For every $\pi \in S_{n+q-1}$, $\Del_q(\pi)-q+1 =
\Del_S(f_q(\pi))$, $\Des_q(\pi)-q+1 = \Des_S(f_q(\pi))$,
$\ell_q(\pi) = \ell_S(f_q(\pi))$, and
$\rmaj_{q,n+q-1}(\pi)=\rmaj_{S_n}(f_q(\pi))$. Here, $X-r = \{\,x-r
\st x\in X\,\}$.
\end{prop}

\begin{prop}[{see~\cite[Proposition~8.4]{rr2}}]\label{PR:fqInv}
For any permutation $w$, $f_q(w)^{-1} = f_q(w^{-1})$.
\end{prop}

The map $f_q$ is obviously not injective for $q>1$. The family of
maps $g_{q,u}$ defined next serve as local inverses of $f_q$ (see
Remark~\ref{REM:gqInvertible}).

\begin{defn}
For $u \in S_{n+q-1}$ with $S$-canonical presentation $u= u_1
\cdots u_{n+q-2}$, define $g_{q,u}:R_j^S \to R_{j+q-1}^S$ by
\[
    g_{q,u}(s_j s_{j-1}\cdots s_\ell) = s_{j+q-1} s_{j+q-2} \cdots s_{\ell+q-1}, \quad
        g_u(s_j s_{j-1}\cdots s_1) = u_{j+q-1}.
\]
Now extend $g_{q,u}:S_n \to S_{n+q-1}$ as follows: let $w\in S_n$,
$w=w_1 \cdots w_{n-1}$ its $S$-canonical presentation, then
\[
    g_{q,u}(w) := u_1 \cdots u_{q-1} \cdot g_{q,u}(w_1)\cdots g_{q,u}(w_{n-1}),
\]
which is clearly the $S$-canonical presentation of $g_{q,u}(w)$.
\end{defn}

\begin{rem}\label{REM:gqInvertible}
Let $w \in S_n$ and $u\in S_{n+q-1}$. Then $f_q(g_{q,u}(w))=w$ if
for all $1 \le j \le n-1$,
\[
    w_j = s_j\cdots s_1 \; \implies \; u_{j+q-1} = s_{j+q-1}\cdots s_\ell,\; \ell \le q   ,
\]
where $w=w_1\cdots w_{n-1}$ and $u=u_1\cdots u_{n+q-2}$ are the
$S$-canonical presentations of $w$ and $u$ respectively.
\end{rem}

\subsection{The map $\Psi_q$}
\begin{defn}
Define $\Psi_q:S_{n+q-1} \to S_{n+q-1}$ by $\Psi_q(v) =
g_{q,v}(\rtlPhi(f_q(v)))$ .
\end{defn}

That is, the image of $v$ under $\Psi_q$ is obtained by applying
$\rtlPhi$ to $f_q(v)$ in $S_n$, then using $g_{q,v}$ as an
``inverse'' of $f_q$ in order to ``lift'' the result back to
$S_{n+q-1}$.

\begin{thm}\label{PR:Psiq}
\begin{enumerate}
\item
The mapping $\Psi_q$ is a bijection of $S_{n+q-1}$ onto itself.
\item
For every $v\in S_{n+q-1}$, $\rmaj_{q,n+q-1}(v) =
\ell_q(\Psi_q(v))$.
\item
For every $v\in S_{n+q-1}$, $\Del_q(v^{-1}) =
\Del_q(\Psi_q(v)^{-1})$.
\item
For every $v\in S_{n+q-1}$,
$\Des_q(v^{-1})=\Des_q(\Psi_q(v)^{-1})$.
\end{enumerate}
\end{thm}

The proof is given below.

\begin{lem}\label{LEM:rep-ltrmq}
Let $v \in S_{n+q-1}$, $v=v_1 \cdots v_{n+q-2}$ its $S$-canonical
presentation. Then for every $q<j \le n+q-1$, $j\in\ltrm_q(v)$ if
and only if $v_{j-1} = s_{j-1} s_{j-2}\cdots s_\ell$ for some
$\ell \le q$.
\end{lem}
\begin{proof}
By induction on $n$. Let $\pi = v_1 \cdots v_{n-1+q-2} \in
S_{n+q-2} \subseteq S_{n+q-1}$ and assume that the assertion is
true for $\pi$. If $v_{n+q-2}=1$, then the claim is correct by the
induction hypothesis. Otherwise,
$v_{n+q-2}=s_{n+q-2}s_{n+q-3}\cdots s_\ell$ for some $1 \le \ell
\le n+q-2$. Writing $\pi=[b_1,\dots,b_{n+q-2}]$, we have that
$v=\pi v_{n+q-2} =
[b_1,\dots,b_{\ell-1},n+q-1,b_\ell,\dots,b_{n+q-2}]$, so clearly
for every $1 \le k \le n+q-2$, the set of numbers smaller than
$b_k$ and to its left in $\pi$ is equal to the set of numbers
smaller than $b_k$ and to its left in $v$. Thus $b_k \in
\ltrm_q(v)$ if and only if $b_k \in \ltrm_q(\pi)$, which, by the
induction hypothesis, is true if and only if $v_{b_k-1} =
s_{b_k-1} \cdots s_r$ for some $r \le q$. Finally, $n+q-1 \in
\ltrm_q(v)$ if and only if $n+q-1$ occupies one of the $q$
leftmost places in $v$, that is, if and only if $\ell \le q$.
\end{proof}

\begin{lem}\label{LEM:fPsiq}
Let $v \in S_{n+q-1}$. Then $f_q(\Psi_q(v)) = \rtlPhi(f_q(v))$.
\end{lem}
\begin{proof}
Let $v = v_1\cdots v_{n+q-2}$ and $w=\rtlPhi(f_q(v))=w_1\cdots
w_{n-1}$ be the $S$-canonical presentations of $v$ and
$\rtlPhi(f_q(v))$ respectively. By definition of $f_q$ and
Corollary~\ref{COR:revPhi}, for every $1 \le j \le n-1$, $w_j =
s_j s_{j-1}\cdots s_1$ if and only if $v_{j+q-1} = s_{j+q-1}\cdots
s_\ell$, $\ell\le q$. Therefore, by Remark~\ref{REM:gqInvertible},
\[
    f_q(\Psi_q(v))=f_q(g_{q,v}(\rtlPhi(f_q(v))))=f_q(g_{q,v}(w))=w = \rtlPhi(f_q(v)) . \qedhere
\]
\end{proof}

\begin{proof}[Proof of Theorem~\ref{PR:Psiq}]
\begin{enumerate}
\item
To prove that $\Psi_q$ is a bijection, it suffices to find its
inverse. Let $v \in S_{n+q-1}$, and let $v=v_1\cdots v_{n+q-2}$,
$w=\rtlPhi(f_q(v))=w_1\cdots w_{n-1}$ and
$u=\Psi_q(v)=g_{q,v}(w)=v_1\cdots v_{q-1} u_q\cdots u_{n+q-2}$ be
the $S$-canonical presentations of $v$, $\rtlPhi(f_q(v))$ and
$\Psi_q(v)$ respectively. By Lemma~\ref{LEM:fPsiq},
\[
\rtlPhi^{-1}(f_q(\Psi_q(v))) =
\rtlPhi^{-1}(\rtlPhi(f_q(v)))=f_q(v),
\]
so
\begin{multline}\label{EQ:inverseq}
    g_{q,\Psi_q(v)}(\rtlPhi^{-1}(f_q(\Psi_q(v)))) = g_{q,\Psi_q(v)}(f_q(v))
  \\= g_{q,u}(f_q(v)) = v_1\cdots v_{q-1} \cdot g_{q,u}(f_q(v_1))\cdots g_{q,u}(f_q(v_{n-1}))  \tag{$*$}     .
\end{multline}
We claim that $\pi \mapsto g_{q,\pi}(\rtlPhi^{-1}(f_q(\pi)))$ is
the inverse of $\Psi_q$, or in other words, that the right hand
side of~\eqref{EQ:inverseq} equals $v_1 v_2 \cdots v_{n+q-2}$. Let
$q \le j \le n+q-2$, and write $v_j = s_j s_{j-1} \cdots s_\ell$.
If $\ell>q$, then $g_{q,u}(f_q(v_j))=g_{q,u}(s_{j-q+1} \cdots
s_{\ell-q+1})=s_j \cdots s_\ell=v_j$. If $\ell \le q$, then
$f_q(v_j)=s_j \cdots s_1$, so by Corollary~\ref{COR:revPhi},
$w_j=s_j \cdots s_1$, and therefore $u_j=g_{q,v}(w_j)=v_j$, so
again $g_{q,u}(f_q(v_j))=v_j$, and the claim is proved.

\item
By Proposition~\ref{PR:fqPairs} and Lemma~\ref{LEM:fPsiq},
$\ell_q(\Psi_q(v)) = \ell_S(f_q(\Psi_q(v))) =
\ell_S(\rtlPhi(f_q(v)))$. By Theorem~\ref{PR:revPhi2} and
Proposition~\ref{PR:fqPairs}, $\ell_S(\rtlPhi(f_q(v))) =
\rmaj_{S_n}(f_q(v)) = \rmaj_{q,n+q-1}(v)$. Thus $\ell_q(\Psi_q(v))
= \rmaj_{q,n+q-1}(v)$ as desired.

\item
By Lemma~\ref{LEM:rep-ltrmq} and the definition of $f_q$,
$\ltrm_q(\Psi_q(v)) = \ltrm(f_q(\Psi_q(v)))-q+1$ (with the
notation $X-r = \{\,x-r \st x\in X\,\}$). Therefore by
Lemmas~\ref{LEM:fPsiq} and~\ref{LEM:ltrm}, $\ltrm_q(\Psi_q(v)) =
\ltrm(\rtlPhi(f_q(v)))-q+1 = \ltrm(f_q(v))-q+1$. Again by
Lemma~\ref{LEM:rep-ltrmq}, we get that $\ltrm_q(\Psi_q(v)) =
\ltrm_q(v)$. By Proposition~\ref{PR:Delk}, this implies that
$\Del_q([\Psi_q(v)]^{-1})\cup\{1,\dots,q\} =
\Del_q(v^{-1})\cup\{1,\dots,q\}$, hence $\Del_q([\Psi_q(v)]^{-1})
= \Del_q(v^{-1})$ as desired.

\item
By Propositions~\ref{PR:fqPairs} and~\ref{PR:fqInv} and
Lemma~\ref{LEM:fPsiq},
\[
\begin{split}
\Des_q([\Psi_q(v)]^{-1})-q+1 &= \Des_S(f_q([\Psi_q(v)]^{-1}))\\
&= \Des_S([f_q(\Psi_q(v))]^{-1})\\
&= \Des_S([\rtlPhi(f_q(v))]^{-1}).
\end{split}
\]
By Remark~\ref{REM:rc}, $[\rtlPhi(f_q(v))]^{-1} = [\rev \Phi
\rev f_q(v)]^{-1} = \com([\Phi \rev f_q(v)]^{-1})$, so
$\Des_S([\rtlPhi(f_q(v))]^{-1}) = \{1,\dots,n-1\} \setminus
\Des_S([\Phi \rev f_q(v)]^{-1})$. By Theorem~\ref{THM:foata},
\[\Des_S([\Phi\rev f_q(v)]^{-1}) = \Des_S([\rev f_q(v)]^{-1}).\]
Hence, $\Des_S([\rtlPhi(f_q(v))]^{-1}) = \{1,\dots,n-1\} \setminus
\Des_S([\rev f_q(v)]^{-1}) = \Des_S(\com([\rev f_q(v)]^{-1}))$.
Since $\com([\rev f_q(v)]^{-1}) = \com(\com([f_q(v)]^{-1})) =
[f_q(v)]^{-1}$, we get that
\[\Des_S([\rtlPhi(f_q(v))]^{-1}) = \Des_S([f_q(v)]^{-1}).\] Finally, by
Propositions~\ref{PR:fqInv} and~\ref{PR:fqPairs},
$\Des_S([f_q(v)]^{-1})=\Des_S(f_q(v^{-1}))=\Des_q(v^{-1})-q+1$.
\qedhere
\end{enumerate}
\end{proof}


\begin{thebibliography}{BW91}
\bibitem[BW91]{bjo}
A. Bj{\"o}rner and Michelle~L. Wachs,
\newblock Permutation statistics and linear extensions of posets.
\newblock {\em J. Combin. Theory Ser. A}, 58(1):85--114, 1991.

\medskip

\bibitem[Car54]{car1}
 L. Carlitz,
\newblock  $q$-Bernoulli and Eulerian numbers,
\newblock {\em Trans. Amer. Math Soc.}, 76 (1954) 332-350.

\medskip

\bibitem[Car75]{car2}
 L. Carlitz,
\newblock  A combinatorial property of $q$-Eulerian numbers,
\newblock {\em Amer. Math. Monthly} 82 (1975) 51-54.

\medskip


\bibitem[Foa68]{foata:netto}
D. Foata,
\newblock On the {N}etto inversion number of a sequence.
\newblock {\em Proc. Amer. Math. Soc.}, 19:236--240, 1968.

\medskip


\bibitem[FH04]{foata:further}
D. Foata and G. N. Han,
\newblock Further properties of the second fundamental transformation on words.
\newblock
  \href{http://www-irma.u-strasbg.fr/~foata/paper/pub92.html}{\texttt{http://w%
ww-irma.u-strasbg.fr/~foata/paper/pub92.html}}, 2004.

\medskip


\bibitem[FS78]{foasch}
D. Foata and M. P. Schutzenberger,
\newblock Major index and inversion number of permutations.
\newblock {\em Math. Nachr. 83 } (1978) 143-159.

\medskip

\bibitem[GG79]{gg}
A. Garsia and I. Gessel,
\newblock Permutation statistics and partitions,
\newblock {\em Adv. Math. 31 } (1979) 288-305.

\medskip

\bibitem[Gol93]{goldschmidt:characters}
David~M. Goldschmidt,
\newblock {\em Group characters, symmetric functions, and the {H}ecke algebra},
  volume~4 of {\em University Lecture Series}.
\newblock American Mathematical Society, Providence, RI, 1993.

\medskip


\bibitem[Kra95]{krat}
C.~Krattenthaler,
\newblock {\em The major counting of non intersecting lattice paths and generating
functions for tableaux}, {\em Mem. Amer. Math.Soc. 115 (1995), no.
552 }.

\medskip


\bibitem[Lot83]{lothaire:words}
M.~Lothaire,
\newblock {\em Combinatorics on words}, volume~17 of {\em Encyclopedia of
  Mathematics and its Applications}.
\newblock Addison-Wesley Publishing Co., Reading, Mass., 1983.

\medskip



\bibitem[MM16]{macmahon16}
P.A.~MacMahon,
\newblock {\em Combinatory Analysis I-II}, {\em Cambridge Univ. Press}.
London/New-York, 1916.
\newblock Reprinted by Chelsea, New-York, 1960.

\medskip

\bibitem[Mit01]{mitsuhashi:alternating}
Hideo Mitsuhashi,
\newblock The {$q$}-analogue of the alternating group and its representations.
\newblock {\em J. Algebra}, 240(2):535--558, 2001.

\medskip

\bibitem[RR03]{rr2}
A. Regev and Y. Roichman,
\newblock $q$ statistics on ${S}_n$ and pattern avoidance.
\newblock
  arXiv:\href{http://arxiv.org/abs/math.CO/0404354}{\texttt{math.CO/0305393}},
  2003.

\medskip

\bibitem[RR04]{regev:alternating}
A. Regev and Y. Roichman,
\newblock Permutation statistics on the alternating group.
\newblock {\em Adv. in Appl. Math.}, 33(4):676--709, 2004.

\medskip

\bibitem[Sta02]{stanley}
R. P. Stanley,
\newblock Some remarks on sign-balanced and maj-balanced posets,
preprint,
\newblock
  arXiv:\href{http://arxiv.org/abs/math.CO/0211113}{\texttt{math.CO/0211113}},
2002.

\end{thebibliography}
\end{document}